\documentclass[11pt,3p,times]{elsarticle}

\usepackage{srcltx}
\usepackage{xcolor}
\usepackage{colortbl}
\usepackage{enumitem}
\usepackage{amssymb, color}
\usepackage{hyperref}
\usepackage{amsmath}
\usepackage{pifont}
\usepackage[utf8]{inputenc}
\makeatletter
\def\LaTeX{\leavevmode L\raise.42ex
    \hbox{\kern-.3em\size{\sf@size}{0pt}\selectfont A}\kern-.15em\TeX}
\makeatother

\newcommand{\BibTeX}{{\rm B\kern-.05em{\sc
          i\kern-.025emb}\kern-.08em\TeX}}
\journal{}
\makeatletter
\def\@currentlabel{2.1}\label{e:dispaa}
\def\@currentlabel{2.21}\label{e:dispau}
\def\@currentlabel{2.22}\label{e:dispav}
\def\@currentlabel{2.23}\label{e:dispaw}
\def\@currentlabel{2.24}\label{e:dispax}
\def\theequation{\thesection.\@arabic\c@equation}
\makeatother

\renewcommand{\theequation}{\arabic{section}.\arabic{equation}}
\newcommand{\R}{\mathbb R}

\def \O{\Omega}
\everymath{\displaystyle}
\newtheorem{theorem}{Theorem}[section]
\newtheorem{thm}{Theorem} [section]
\newtheorem{lem}{Lemma} [section]
\newtheorem{prop}{Proposition} [section]

\newtheorem{definition}{Definition} [section]

\newtheorem{rem}{Remark}[section]

\renewcommand{\theequation}{\thesection.\arabic{equation}}
\renewcommand{\thesection}{\arabic{section}}
\renewcommand{\theequation}{\thesection.\arabic{equation}}
\let\ssection=\section\renewcommand{\section}{\setcounter{equation}{0}\ssection}
\begin{document}
\begin{frontmatter}
\title{A new class of multiple nonlocal problems with two parameters and variable-order fractional $p(\cdot)$-Laplacian}

\author[mk0,mk1,mk2]{Mohamed Karim Hamdani\corref{cor1}}
\ead{hamdanikarim42@gmail.com}
\author[L]{Lamine Mbarki}
\ead{mbarki.lamine2016@gmail.com}
\author[OMA,ALL]{Mostafa Allaoui}
\ead{m.allaoui@uae.ac.ma}
\begin{center}
\address[mk0]{Science and technology for defense lab LR19DN01, center for military research, military academy, Tunis, Tunisia.}
\address[mk1]{Military Aeronautical Specialities School, Sfax, Tunisia.}
\address[mk2]{Department of Mathematics, University of Sfax, Faculty of Science of Sfax, Sfax, Tunisia.}
\address[L]{Mathematics Departement, Faculty of Science of Tunis , University of Tunis El Manar, Tunisia.}
\address[OMA]{Department of Mathematics, FSTH Abdelmalek Essaadi University-Tetuan, Morocco.}
\address[ALL]{Department of Mathematics, Mohammed I University,  Oujda, Morocco.}
\end{center}
\begin{abstract}
In the present manuscript, we focus on a novel tri-nonlocal Kirchhoff problem, which involves the $p(x)$-fractional Laplacian equations of variable order. The problem is stated as follows:
\begin{eqnarray*}
\left\{
 \begin{array}{ll}
 M\Big(\sigma_{p(x,y)}(u)\Big)(-\Delta)^{s(\cdot)}_{p(\cdot)}u(x)
=\lambda |u|^{q(x)-2}u\left(\int_\O\frac{1}{q(x)} |u|^{q(x)}dx \right)^{k_1}+\beta|u|^{r(x)-2}u\left(\int_\O\frac{1}{r(x)} |u|^{r(x)}dx \right)^{k_2}
 \quad \mbox{in }\Omega, \\
 \\
u=0     \quad \mbox{on }\partial\Omega,
\end{array}
\right.
\end{eqnarray*}
where the nonlocal term is defined as
$$
\sigma_{p(x,y)}(u)=\int_{\Omega\times \Omega}\frac{1}{p(x,y)}\frac{|u(x)-u(y)|^{p(x,y)}}{|x-y|^{N+s(x,y)p(x,y)}} \,dx\,dy.
$$
Here, $\Omega\subset\mathbb{R}^{N}$ represents a bounded smooth domain with at least $N\geq2$. The function $M(s)$ is given by $M(s) = a - bs^\gamma$, where $a\geq 0$, $b>0$, and $\gamma>0$. The parameters $k_1$, $k_2$, $\lambda$ and $\beta$ are real parameters, while the variables $p(x)$, $s(\cdot)$, $q(x)$, and $r(x)$ are continuous and can change with respect to $x$.
To tackle this problem, we employ some new methods and variational approaches along with two specific methods, namely the Fountain theorem and the symmetric Mountain Pass theorem. By utilizing these techniques, we establish the existence and multiplicity of solutions for this problem separately in two distinct cases: when $a>0$ and when $a=0$. To the best of our knowledge, these results are the first contributions to research on  the variable-order $p(x)$-fractional Laplacian operator.
\end{abstract}

\end{frontmatter}

\section{Statement of the Problem and the Main Results}\label{sect1}
Given that $ N\geq 2 $ and $\Omega\subset  \mathbb{R}^{N}$, is a smooth bounded domain. The goal of this paper is to investigate the existence and multiplicity of solutions for variable order $p(x)$-Kirchhoff tri-nonlocal fractional equations.
 \begin{eqnarray}\label{main}
\left\{
 \begin{array}{ll}
 M\Big(\sigma_{p(x,y)}(u)\Big)(-\Delta)^{s(\cdot)}_{p(\cdot)}u(x)
=\lambda |u|^{q(x)-2}u\left(\int_\O\frac{1}{q(x)} |u|^{q(x)}dx \right)^{k_1}+\beta|u|^{r(x)-2}u\left(\int_\O\frac{1}{r(x)} |u|^{r(x)}dx \right)^{k_2}
 \quad \mbox{in }\Omega, \\
 \\
u=0     \quad \mbox{on }\partial\Omega,
\end{array}
\right.
\end{eqnarray}
where
$$
\sigma_{p(x,y)}(u)=\int_{\Omega\times \Omega}\frac{1}{p(x,y)}\frac{|u(x)-u(y)|^{p(x,y)}}{|x-y|^{N+s(x,y){p(x,y)}}} \,dx\,dy,
$$
where   $N > s(x,y)p(x, y)$ for all $(x, y)\in \overline{\Omega}\times \overline{\Omega}$, $\lambda,\beta$ are two real parameters, $k_1,k_2>0$, $M(x)=a-bx^{\gamma}$, $a\geq 0$, $b, \gamma>0$ and $q$, $r$ are  continuous real functions on $\bar{\Omega}$.\\
The operator defined as $(-\Delta)^{s(\cdot)}_{p(\cdot)}$ is referred to as the $p(x)$-fractional Laplacian with variable order, and it is defined as follows:
$$(-\Delta)^{s(\cdot)}_{p(\cdot)}u(x):=P.V.\int_{\Omega}\frac{1}{p(x,y)}\frac{|u(x)-u(y)|^{p(x,y)-2}(u(x)-u(y))}{|x-y|^{N+s(x,y){p(x,y)}}}\,dy;$$
for any $u\in C_{0}^{\infty} (\mathbb{R}^N)$, where the notation P.V. means the Cauchy principal value.

As the problem \eqref{main} involves integrals over the domain $\Omega$, it deviates from being a pointwise identity. Consequently, it is commonly referred to as a tri-nonlocal problem due to the presence of the following integrals.
$$\sigma_{p(x,y)}(u) \mbox{ and }\int_\O\frac{1}{s(x)} |u|^{s(x)}dx,\;\mbox{for } s=\{q,r\}.$$
In recent years, the wide class of problems involving nonlocal operators have been an increasing attention and have acquired a renovate relevance due to their occurrence in pure and applied mathematical point view, for instance, the finance, thin obstacle problem, biology as the interaction of bacteria, probability, optimization and others.

In the current work, our attention will be focused on a very interesting nonlocal operator known as the fractional $p(x)$-Laplacian with variable order. This type of operator represents an extension and a combination of many other operators. Indeed, the nonlocal fractional $p$-Laplacian, which has been extensively studied in the literature, is defined as
\[
(-\Delta)_p^s u(x) = 2\lim_{\varepsilon \to 0^+} \int_{\mathbb{R}^N \setminus B_\varepsilon(0)} \frac{|u(x) - u(y)|^{p-2}(u(x) - u(y))}{|x - y|^{N+sp}} dy, \quad x \in \mathbb{R}^N.
\]
During this time, problems involving variable exponents have attracted many researchers \cite{DHHR, FZ2001, GTW}. These types of problems primarily arise from the $p(x)$-Laplace operator $\operatorname{div}(|\nabla u|^{p(x)-2}\nabla u)$, which serves as a natural extension of the classical $p$-Laplace operator $\operatorname{div}(|\nabla u|^{p-2}\nabla u)$ when $p$ is a positive constant. However, these operators possess a more intricate structure due to their lack of homogeneity. Hence, problems involving $p(x)$-Laplacian become more tricky. Moreover, concerning the nonlocal problem involving the $p(x)$-Laplacian, we can refer to \cite{JSC,HHMR,HCR,HLAD, JAB, ENR, EAO, AMHB} and the references therein. For instance, in \cite{JAB}, the authors focused their study on a specific fourth-order bi-nonlocal elliptic equation of Kirchhoff type with Navier boundary conditions, which is expressed as:

\begin{eqnarray*}
\left\{
 \begin{array}{ll}
 M\Big(\int_{\Omega}\frac{1}{p(x)|\Delta u|^{p(x)dx}}\Big)\Delta_{p(x)}^{2}u(x)
=\lambda |u|^{q(x)-2}u\left(\int_\O\frac{1}{q(x)} |u|^{q(x)}dx \right)^{r}\quad \mbox{in }\Omega, \\
 \\
\Delta u=u=0     \quad \mbox{on }\partial\Omega,
\end{array}
\right.
\end{eqnarray*}
By using a variational method and critical point theory, the authors obtained a nontrivial weak solution. Consequently, the idea to replace the fractional $p$-Laplacian by its variable version was initiated. For this purpose, Kaufmann et al. \cite{KRV} introduced the fractional $p(x)$-Laplacian $(-\Delta)_{p(\cdot)}^{s}$ as follows:
\[
(-\Delta)_{p}^{s}u(x) = \lim_{\varepsilon \to 0^{+}}\int_{\mathbb{R}^{N}\setminus B_{\varepsilon}(0)}\frac{|u(x)-u(y)|^{p(x,y)-2}(u(x)-u(y))}{|x-y|^{N+sp(x,y)}}dy,\quad x\in \mathbb{R}^{N}.
\]
To address such problems, the authors considered the fractional Sobolev space with variable exponents, variational methods, existence. Simultaneously, many works involving the variable-order fractional Laplacian (see \cite{XZY}) have emerged, defined as follows:
\[
(-\Delta)^{s(\cdot)}u(x) = \lim_{\varepsilon \to 0^{+}}\int_{\mathbb{R}^{N}\setminus B_{\varepsilon}(0)}\frac{|u(x)-u(y)|}{|x-y|^{N+2s(x,y)}}dy,\quad x\in \mathbb{R}^{N}.
\]
Furthermore, the combination of these operators leads to the emergence of the so-called fractional $p(x)$-Laplacian with variable order. This class of operators has captured the attention of numerous researchers \cite{AHL, GY, WQHKY, XZY, ZYL, ZFB}, who have investigated various aspects, including the existence, multiplicity, and qualitative properties of the solutions. Additionally, there are several works focusing on the nonlocal fractional $p(x)$-Laplacian with variable order \cite{AHL, BT(EJDE), BT(TMNA2), GY, XZY, ZYL} and their references therein. For instance, in \cite{ZYL}, the authors studied the existence and multiplicity of solutions for the following fractional $p(\cdot)$-Kirchhoff type problem with a variable order $s(\cdot)$:
\begin{equation}\label{eqrr}
\left\{
\begin{array}{ll}
M\left(\iint_{\mathbb{R}^{2N}}\frac{1}{p(x,y)}\frac{|v(x)-v(y)|^{p(x,y)}}{|x-y|^{N+p(x,y)s(x,y)}}dxdy\right)(-\Delta)^{s(\cdot)}_{p(\cdot)}v(x)+|v(x)|^{\overline{p}(x)-2}v(x) =\mu g(x,v) & \text{in }\mathbb{R}^{N},\\
v\in W^{s(\cdot),p(\cdot)}(\mathbb{R}^{N}),
\end{array}\right.
\end{equation}
where  $(x,y)\in\mathbb{R}^{N}\times\mathbb{R}^{N}$ satisfies the condition $N>p(x,y)s(x,y)$, $s(\cdot):\mathbb{R}^{2N}\to(0,1)$ and $p(\cdot):\mathbb{R}^{2N}\to(1,\infty)$, and $\overline{p}(x)=p(x,x)$ for $x\in\mathbb{R}^{N}$, $M$ is a continuous Kirchhoff-type function, $g(x,v)$ is a Carath\'{e}odory function and $\mu>0$ is a parameter. The authors obtained at least two distinct solutions for the above problem by applying the generalized abstract critical point theorem. In addition, under weaker conditions, they also proved the existence of one solution and infinitely many solutions using the mountain pass lemma and fountain theorem, respectively.

Motivated by the aforementioned works, the present work aims to study the problem \eqref{main} mentioned above. The main difficulties and innovations lie in the form of the new Kirchhoff functions $M(s) = a - bs^\gamma$, derived from the negative Young's modulus when the atoms are spread apart rather than compressed together, resulting in negative deformation. In the case $a=0$, to overcome this challenge, inspired by \cite{ABMS}, our main approach is based on the notion of the first eigenvalue associated with our operator.

The specificity of this tool is that in the literature, we only find the recent paper \cite{ABMS}, in which the authors introduce the $s(\cdot,\cdot)$-fractional Musielak-Sobolev spaces $W^{s(x,y)}L_{\varphi(x,y)}(\Omega)$. By employing Ekeland's variational principle, the authors establish the existence of a positive value $\lambda^{**}>0$ such that for any $\lambda$ within the interval $(0,\lambda^{**})$, it serves as an eigenvalue for the following problem:

\[
(\mathcal{P}_a) \quad \left\{
\begin{array}{clclc}
\left(-\Delta\right)^{s(x,.)}_{a_{(x,.)}} u & = & \lambda |u|^{q(x)-2}u & \text{in }\Omega, \\
u & = & 0 & \text{in }\mathbb{R}^N\setminus \Omega,
\end{array}
\right.
\]
where $\Omega$ is a bounded open subset of $\mathbb{R}^N$ with a $C^{0,1}$-regularity and a bounded boundary conditions. It is noteworthy that this operator represents a generalization of $(-\Delta)^{s(\cdot)}_{p(\cdot)}$ (whenever we take $a_{(x,.)}=t^{p(x,.)-2}$). Thus, this characterization is applicable in our case.

Throughout this work, the functions $s(\cdot)$ and $p(\cdot)$ satisfy the following conditions:

\begin{enumerate}[label={$(H\sb{\arabic*})$:},ref={{$(H\sb{\arabic*})$}}]
\item\label{hypot:S} $s(x,y)$ is a symmetric function, i.e., $s(x,y)=s(y,x)$, and we have
  \[0 < s^{-} := \inf_{(x,y)\in \overline{\Omega}\times \overline{\Omega}} s(x,y) \leq s^{+} := \sup_{(x,y)\in \overline{\Omega}\times \overline{\Omega}} s(x,y) < 1.\]

  \item\label{hypot:P} $p(x,y)$ is a symmetric function, i.e., $p(x,y)=p(y,x)$, and we have
  \[1 < p^{-} := \inf_{(x,y)\in \overline{\Omega}\times \overline{\Omega}} p(x,y) \leq p^{+} := \sup_{(x,y)\in \overline{\Omega}\times \overline{\Omega}} p(x,y) < \infty.\]
\end{enumerate}

For any $x\in\Omega$, we denote
\[\overline{p}(x) := p(x,x), \quad \overline{s}(x) := s(x,x).\]

Moreover, it is also assumed that the function $p(\cdot)$ satisfies the following condition:
\begin{equation}\label{cond-pqr}
  \min\left\{ p^-(\gamma+1), q^-(k_1+1), r^-(k_2+1) \right\} > p^+.
\end{equation}
Now, we are ready to state our main results.
\begin{thm}\label{theo1.1}
Let $a > 0$. Assume that \ref{hypot:S}--\ref{hypot:P} and \eqref{cond-pqr} are satisfied. Then, for all $\lambda > 0$ and $\beta > 0$, the problem \eqref{main} possesses a nontrivial weak solution in $X$.
\end{thm}

\begin{thm}\label{theo1.2}
Let $a > 0$. Assume that \ref{hypot:S}--\ref{hypot:P} and \eqref{cond-pqr} are satisfied. Then, for all $\lambda > 0$ and $\beta > 0$, the problem \eqref{main} admits an unbounded sequence of solutions in $X$.
\end{thm}

Moreover, we also obtain the following existence results for problem \eqref{main} in the case $a = 0$.

\begin{thm}\label{symm.moun.pass}
Let $a = 0$, $\lambda_1$ be the first eigenvalue associated with our operator (For more comprehensive details and references, we recommend referring to \cite{ABMS}), and
\begin{align*}
\left\{
 \begin{array}{ll}
   \lambda^* = -\frac{\beta\lambda_1(x)(q^-)^{k_1+1}(k_1+1)}{4C(r^+)^{k_2+1}(k_2+1)\rho^{q^-(k_1+1)-r^+(k_2+1)}}, & \text{if } \lambda > 0, \beta < 0,\\\\
   \beta^* = -\frac{\lambda\lambda_1(x)(r^-)^{k_2+1}(k_2+1)}{4C(q^+)^{k_1+1}(k_1+1)\rho^{r^-(k_2+1)-q^+(k_1+1)}}, & \text{if } \lambda < 0, \beta > 0.
 \end{array}
\right.
\end{align*}
If any of the following conditions is satisfied:
\begin{align}\label{cond(1.4)a=0}
\lambda < 0, \beta \in (0, \beta^*) \text{ and } r^+(k_2+1) < \min\left\{ (\gamma + 1)p^-, q^-(k_1+1) \right\},\nonumber\\
\lambda \in (0, \lambda^*), \beta < 0 \text{ and } q^+(k_1+1) < \min\left\{ (\gamma + 1)p^-, r^-(k_2+1) \right\},\\
\lambda < 0, \beta < 0 \text{ and } r^+(k_2+1) < q^-(k_1+1) < (\gamma + 1)p^-,\nonumber
\end{align}
then problem \eqref{main} has infinitely many solutions in $X$.
\end{thm}
However, as far as our knowledge extends, there are no existing results regarding the existence and multiplicity of solutions for problem \eqref{main} involving the new tri-nonlocal Kirchhoff function and the $p(x)$-fractional Laplacian operator with variable order.

The structure of this paper is as follows: In the second section, an abstract framework is presented, where we provide a review of some preliminary results that will be utilized throughout the subsequent sections. The third section is specifically focused on presenting the Palais-Smale condition separately for the cases of $a > 0$ and $a = 0$. The subsequent sections are dedicated to proving the main results of this study.
\section{Abstract framework}
\subsection{\textbf{Generalized Lebesgue and Sobolev spaces}}
In this section, we provide a brief review of the definition and key results concerning Lebesgue spaces with variable exponents and generalized Sobolev spaces. For a more comprehensive understanding, interested readers are referred to \cite{DHHR, FZ2001, R} and the references therein.

For this purpose, let us define
\[
C_+(\Omega) := \{ h : h \in C(\overline{\Omega}) \text{ and } h(x) > 1 \text{ for all } x \in \overline{\Omega} \}.
\]
For $p(\cdot)\in C_+(\Omega)$,  the variable exponent Lebesgue
space $L^{p(\cdot)}(\Omega)$ is defined by
$$
L^{p(\cdot)}(\Omega):=\{u:\Omega\to\mathbb{R}\text{ measurable and }
 \int_{\Omega}|u(x)|^{p(x)}dx<\infty\}.
$$
This space is endowed with the so-called Luxemburg norm given by
$$
\| u \|_{L^{p(x)}(\Omega)}=|u|_{p(\cdot)}:=\inf\{\delta>0:\int_{\Omega}|\frac{u(x)}{\delta}|^{p(x)}dx\leq1\}
$$
and $( L^{p(\cdot)}(\Omega),|u|_{p(\cdot)})$ becomes a Banach space, and we call it variable exponent Lebesgue space.

Now, in order to claim the $(PS)$ condition cited in Section\ref{sec3}, we  state the following lemma for the variable exponent Lebesgue spaces (see \cite[Lemma A.1]{GTW})

\begin{lem}\label{lemma2.2}
Assume that $h_{1}\in L^{\infty}(\Omega)$ such that $h_{1}\geq0$ and $h_{1}\not\equiv0$ a.e. in $\Omega$. Let $h_{2}:\Omega\rightarrow\mathbb{R}$ be a measurable function such that $h_{1}h_{2}\geq1$ a.e. in $\Omega$. Then for any $u\in L^{h_{1}(\cdot)h_{2}(\cdot)}(\Omega)$,
\begin{equation*}
\||u|^{h_{1}(\cdot)}\|_{h_{2}(\cdot)}\leq\|u\|^{h_1^-}_{h_1(\cdot)h_2(\cdot)}+\|u\|^{h_1^+}_{h_1(\cdot)h_2(\cdot)}.
\end{equation*}
\end{lem}
The generalized Sobolev space, denoted by
$W^{k,p(\cdot)}(\Omega)$, is defined as follows
$$
W^{k,p(\cdot)}(\Omega)=\{u\in L^{p(\cdot)}(\Omega)|\;D^\alpha
u\in L^{p(\cdot)}(\Omega),|\alpha|\leq k\}
$$
where
$$
D^\alpha u=\frac{\partial^{|\alpha|}}{\partial x_1^{\alpha_1}
\cdots\partial x_N^{\alpha_N}}u
$$
with $\alpha= (\alpha_1,\ldots ,\alpha_N)$ is a multi-index and
$|\alpha| = \sum_{i=1}^{N} \alpha_i$. The space $W^{k,p(\cdot)}(\Omega)$,
equipped with the norm
$$
\|u\|_{k,p(\cdot)}:=\sum_{|\alpha|\leq k}|D^\alpha u|_{p(\cdot)},
$$
is a uniformly convex, separable, and reflexive Banach space.
\subsection{\textbf{Fractional Sobolev spaces with variable exponents}}

In the present part, we recall some properties of the fractional Sobolev spaces with  variable
 exponents which will be useful in the rest of the paper. For more details, we can refer to \cite {BR,BT(CVEE),BT(EJDE),BT(TMNA2),KRV}.

In the present part, we give the variational setting of problem (1.1) and state important results to be used later.   We set $\mathcal{Q}:=\mathbb R^{2N}\setminus(C^\Omega_{R^{N}}\times C^\Omega_{R^{N}})$ and define the
 fractional Sobolev space with variable exponent as
\begin{align*}
 {X}&=W^{s(x,y),\overline{p}(x),p(x,y)}(\Omega)\\
 & :=\Big\{u:\mathbb R^N \to\mathbb R:u_{|_\Omega}\in L^{\overline{p}(x)}(\Omega),\\
 &\quad \int_{\mathcal{Q}}\frac{|u(x)-u(y)|^{p(x,y)}}{\eta^{p(x,y)}| x-y|^{N+s(x,y)p(x,y)}}
 \,dx\,dy<\infty,\text{ for some } \eta>0\Big\}.
 \end{align*}
The space ${X}$ is equipped with the  norm
$$\| u\|_{ {X}}:=\| u\|_{L^{\overline{p}(x)}(\Omega)}
 +[u]_{X};$$

where $[u]_{X}$ is the seminorm defined as follows
$$
[u]_{X}=\inf\Big\{\eta>0:\int_{\mathcal{Q}}\frac{|
 u(x)-u(y)|^{p(x,y)}}{\eta^{p(x,y)}| x-y |^{N+s(x,y)p(x,y)}}\,dx\,dy<1\Big\}.
$$
Then $(X,\|\cdot\|_{X})$ is a separable reflexive Banach space.\\
 Now, let define the subspace ${X}_0$ of ${ X}$ as
$$
 X_0 = :=\{u\in { X}: u=0\text{ a.e.\ in }\Omega^c\}.
$$
We define the norm on ${ X}_0$ as follows
\[
 \| u\|_{{ X}_0}:=\inf\Big\{\eta>0:\int_{\mathcal{Q}}\frac{|
 u(x)-u(y)|^{p(x,y)}}{\eta^{p(x,y)}| x-y |^{N+s(x,y)p(x,y)}}\,dx\,dy<1\Big\}.
\]

 \begin{rem}\label{rem1}
 For $u\in X_0$, we obtain
 $$
\int_{\mathcal{Q}}\frac{|
 u(x)-u(y)|^{p(x,y)}}{\eta^{p(x,y)}| x-y |^{N+s(x,y)p(x,y)}}\,dx\,dy
= \int_{\mathbb R^N\times\mathbb R^N}\frac{| u(x)-u(y)|^{p(x,y)}}
 {\eta^{p(x,y)}| x-y |^{N+s(x,y)p(x,y)}}\,dx\,dy.
$$
Thus, we have
\[
 \| u\|_{{ X}_0}:=\inf\Big\{\eta>0:\int_{\mathbb R^N\times\mathbb R^N}\frac{|
 u(x)-u(y)|^{p(x,y)}}{\eta^{p(x,y)}| x-y |^{N+s(x,y)p(x,y)}}\,dx\,dy<1\Big\}.
\]
\end{rem}

Now we state the following continuous and compact embedding result for the
space ${ X}_0$. The proof follows from \cite[Theorem 2.2, Remark 2.2]{ABS}.

\begin{thm}\label{prp 3.3}
Let $\Omega $ be a smooth bounded domain in $\mathbb{R}^N $, $s(\cdot,\cdot)\in(0,1)$ and $p(\cdot, \cdot)$  satisfy \ref{hypot:S} and \ref{hypot:P} with $s^+p^+<N$.
 Then for any  $r\in C_+(\overline{\Omega})$ such that $1<r(x)< p_s^*(x)$ for all
$x\in \overline{\Omega}$, there exits a constant $C=C(N,s,p,r,\Omega)>0$
such that for every $u\in{ X}_0$,
\[
 \| u \|_{L^{r(x)}(\Omega)}\leq C \| u \|_{{ X}_0}.
\]
Moreover, this embedding is compact.
\end{thm}

\begin{definition} \rm
 For $u\in { X}_0$, we define the  modular  $\rho_{{ X}_0}:{ X}_0\to  \mathbb R$ as
 \begin{equation} \label{modular}
 \rho_{{ X}_0}(u):=\int_{\mathbb R^N \times \mathbb R^N}\frac{|
 u(x)-u(y)|^{p(x,y)}}{| x-y |^{N+s(x,y)p(x,y)}}\,dx\,dy.
\end{equation}
\end{definition}

The interplay between the norm in ${ X}_0$ and the modular function $\rho_{{ X}_0}$ can be studied in the following lemma.

\begin{lem}\label{lem 3.1}
 Let $u \in { X}_0$ and $\rho_{{ X}_0}$ be defined as in \eqref{modular}.
 Then we have the following results:
 \begin{itemize}
 \item[(i)] $ \| u \|_{{ X}_0}<1$ $(=1;>1)$  if and only if $\rho_{{ X}_0}(u)<1(=1;>1)$.

 \item[(ii)] If $\| u \|_{{ X}_0}>1$, then
 $ \| u \|_{{ X}_0} ^{p^{-}}\leq\rho_{{ X}_0}(u)\leq\| u \|_{{ X}_0}^{p{+}}$.

 \item[(iii)] If $\| u \|_{{ X}_0}<1$, then
 $\| u \|_{{ X}_0} ^{p^{+}}\leq\rho_{{ X}_0}(u)\leq\| u \|_{{ X}_0}^{p{-}}$.
 \end{itemize}
\end{lem}	

The next lemma can easily be obtained using the properties of the modular function
 $\rho_{X_0}$  from  Lemma \ref{lem 3.1}.
	
\begin{prop}[{\cite{ZYL, ZAF}}]\label{lem 3.2}
 Let $u,u_m  \in { X}_0$, $m\in\mathbb N$. Then the following two statements are equivalent:
 \begin{itemize}
 \item[(i)] ${\lim_{m\to  \infty} }\| u_m - u \|_{{ X}_0} =0$,
 \item[(ii)] ${\lim_{m\to  \infty}} \rho_{{ X}_0}(u_m -u)=0$.
 \end{itemize}
\end{prop}

\begin{lem}[{\cite[Lemma 2.3]{ABS}}] \label{lem 3.3}
 $(X_0,\|\cdot\|_{X_0})$ is a separable, reflexive and uniformly convex Banach space.
\end{lem}
\section{\textbf{Checking of the $(PS)_c$ condition:}}\label{sec3}
In this part, we will use as space of work, the space $X_0$ and by simplicity we will denote this as $X$ instead of $X_0$ in the rest of this paper.
\\
Considering the variational structure of \eqref{main}, we look for critical points of the corresponding Euler-Lagrange functional $\mathcal{I}_{\lambda, \beta}: X \rightarrow \mathbb{R}$, which is defined as follows:
\begin{eqnarray}\label{e3.1}
\mathcal{I}_{\lambda, \beta}(u) = a\sigma_{p(x,y)}(u)-\frac{b}{\gamma+1}\left(\sigma_{p(x,y)}(u)\right)^{\gamma+1} -\frac{\lambda}{k_1+1}\left(\int_\Omega\frac{1}{q(x)} |u|^{q(x)}dx \right)^{k_1+1}-\frac{\beta}{k_2+1}\left(\int_\Omega\frac{1}{r(x)} |u|^{r(x)}dx \right)^{k_2+1},
\end{eqnarray}

for all $u\in X$. It is important to note that $\mathcal{I}_{\lambda, \beta}$ is a $C^1(X, \mathbb{R})$ functional, and its derivative can be computed as follows:
\begin{eqnarray} \label{deriv-fun}
\langle \mathcal{I}'_{\lambda,\beta}(u),\phi\rangle &=& \left[a-b\left(\sigma_{p(x,y)}(u)\right)^{\gamma}\right]\int_{\Omega \times \Omega}\frac{|u(x)-u(y)|^{p(x,y)-2}(u(x)-u(y))(\phi(x)-\phi(y))}{|x-y|^{N+p(x,y)s(x,y)}}dxdy \nonumber\\&\quad-&\lambda\left(\int_{\O}\frac{1}{q(x)} |u|^{q(x)}dx \right)^{k_1}\int_{\O} |u|^{q(x)-2}u \phi dx - \beta\left(\int_{\O}\frac{1}{r(x)} |u|^{r(x)}dx \right)^{k_2}\int_{\O} |u|^{r(x)-2}u \phi dx\nonumber,\\
\end{eqnarray}
for any $v\in X$. Consequently, critical points of $\mathcal{I}_{\lambda, \beta}$ correspond to weak solutions of \eqref{main}.
\subsection{\textbf{The $(PS)_c$ condition for $\mathbf{a>0}$}}
\begin{lem}\label{lem2}
Assuming that \eqref{cond-pqr} is valid, then the functional $\mathcal{I}_{\lambda, \beta}$ satisfies the Palais-Smale condition at level $c$, where $c \in\left(0, \frac{a^\frac{\gamma+1}{\gamma}}{ b^\frac{1}{\gamma}+\frac{b^\frac{1}{\gamma}}{\gamma}}\right)$.
\end{lem}
{\bf Proof.} Let ${u_n}$ be a $(PS)c$ sequence of $\mathcal{I}{\lambda, \beta}$ with $c \in\left(0, \frac{a^\frac{\gamma+1}{\gamma}}{ b^\frac{1}{\gamma}+\frac{b^\frac{1}{\gamma}}{\gamma}}\right)$. This implies that the following conditions hold:
\begin{eqnarray}\label{er1}
\mathcal{I}_{\lambda, \beta}(u_n) \to c, \quad \mathcal{I}'_{\lambda, \beta}(u_n) \to 0 \mbox{ in } X^*, \quad n \to \infty,
\end{eqnarray}
where $X^*$ denotes the dual space of $X$.
\\
{\bf Step 1.} Firstly, we aim to prove that the sequence ${u_n}$ is bounded in $X$. By assuming the contrary, i.e., supposing that ${u_n}$ is unbounded in $X$, so up to a subsequence, we may assume that $\|u_n\|_X \to \infty$ as $n\to \infty$.  we have
\begin{eqnarray}\label{er2}
& &p^+c+1+\|u_n\|_X  \geq  p^+\mathcal{I}_{\lambda, \beta}(u_n) - \langle \mathcal{I}_{\lambda, \beta}'(u_n),u_n\rangle \notag\\
& = & p^+\left[a\sigma_{p(x,y)}(u_n)-\frac{b}{\gamma+1}\left(\sigma_{p(x,y)}(u_n)\right)^{\gamma+1} -\frac{\lambda}{k_1+1}\left(\int_\Omega\frac{1}{q(x)} |u_n|^{q(x)}dx \right)^{k_1+1}-\frac{\beta}{k_2+1}\left(\int_\Omega\frac{1}{r(x)} |u_n|^{r(x)}dx \right)^{k_2+1}\right]\notag\\
&\qquad -&\Bigg(\left[a-b\left(\sigma_{p(x,y)}(u_n)\right)^{\gamma}\right]\int_{\Omega\times \Omega}\frac{|u_n(x)-u_n(y)|^{p(x,y)}}{|x-y|^{N+p(x,y)s(x,y)}}dxdy -\lambda\left(\int_{\Omega}\frac{1}{q(x)} |u_n|^{q(x)}dx \right)^{k_1}\int_{\Omega} |u_n|^{q(x)}dx\notag\\ & \qquad-& \beta\left(\int_{\Omega}\frac{1}{r(x)} |u_n|^{r(x)}dx \right)^{k_2}\int_{\Omega} |u_n|^{r(x)}dx\Bigg)\notag\\
&\geq&\frac{b}{(p_+)^\gamma}\left(1-\frac{p^+}{(\gamma+1)(p^-)}\right)\left(\int_{\Omega\times \Omega}\frac{|u_n(x)-u_n(y)|^{p(x,y)}}{|x-y|^{N+p(x,y)s(x,y)}}dxdy\right)^{\gamma+1}\notag\\
&\qquad+&\frac{\lambda}{(q^+)^{k_{1}}}\left(1-\frac{p^+}{(k_1+1)(q^-)}\right)\left(\int_{\Omega} |u_n|^{q(x)}dx \right)^{k_1+1}+\frac{\beta}{(r^+)^{k_2}}\left(1-\frac{p^+}{(k_2+1)(r^-)}\right)\left(\int_{\Omega} |u_n|^{r(x)}dx \right)^{k_2+1}.
\end{eqnarray}
From \eqref{cond-pqr} and the fact that $\gamma>0$ and $k_i>0$ for $i=1,2$, it follows that
\begin{eqnarray}\label{er3}
\left\{
 \begin{array}{ll}
1-\frac{p^+}{(\gamma+1)(p^-)}>0, \\
1-\frac{p^+}{(k_1+1)(q^-)}>0,\\
1-\frac{p^+}{(k_2+1)(r^-)}>0.
\end{array}
\right.
\end{eqnarray}
We deduce from \eqref{er2} and \eqref{er3}, that
\begin{gather*}\label{er4}
p^+c+1+\|u_n\|_X \geq \frac{b}{(p^+)^\gamma}\left(1-\frac{p^+}{(\gamma+1)(p^-)}\right)
\left(\int_{\Omega\times \Omega}\frac{|u_n(x)-u_n(y)|^{p(x,y)}}{|x-y|^{N+p(x,y)s(x,y)}}dxdy\right)^{\gamma+1}.
\end{gather*}
If the sequence $(u_n)$ is unbounded in $X$, we can assume, by passing to a subsequence if necessary, that $\|u_n\|_X> 1$. Considering the previous inequalities, we have the following:
\begin{gather*}\label{er5}
p^+c+1+\|u_n\|_X \geq \frac{b}{(p_+)^\gamma}\left(1-\frac{p^+}{(\gamma+1)(p^-)}\right)
\|u_n\|^{(\gamma+1)p^-}_X,
\end{gather*}
which is absurd since $(\gamma+1)p^->1$. Thus, $\{u_n\}$ must be bounded in $X$, and the first assertion is proven.
\\
{\bf Step 2.} Now, we aim to demonstrate that the sequence $\{u_n\}$ has a convergent subsequence in $X$. According to Theorem \ref{prp 3.3}, the embedding $X \hookrightarrow L^{\tau(x)}(\Omega)$ is compact, where $1\leq \tau(x)<p^*_{s}(x)$. Since $X$ is a reflexive Banach space, passing, if necessary, to a subsequence, there exists $u\in X$ satisfying:
\begin{equation}
\label{cvg}
u_n\rightharpoonup u\mbox{ in } X,\; u_n \to u \mbox{ in } L^{\tau(x)}(\Omega),\;\; u_n(x)\to u(x), \mbox{ a.e. in } \Omega.
\end{equation}
From \eqref{deriv-fun}, we find that
\begin{eqnarray} \label{er6}
& &\langle \mathcal{I}_{\lambda, \beta}'(u),u_n-u\rangle\nonumber\\ &=& \left[a-b\left(\sigma_{p(x,y)}(u_n)\right)^{\gamma}\right]\int_{\Omega\times \Omega}\frac{|u_n(x)-u_n(y)|^{p(x,y)-2}(u_n(x)-u_n(y))((u_n(x)-u(x))-(u_n(y)-u(y)))}{|x-y|^{N+p(x,y)s(x,y)}}dxdy \nonumber\\&\quad-&\lambda\left(\int_{\Omega}\frac{1}{q(x)} |u_n|^{q(x)}dx \right)^{k_1}\int_{\Omega} |u_n|^{q(x)-2}u_n (u_n-u) dx -\beta \left(\int_{\Omega}\frac{1}{r(x)} |u_n|^{r(x)}dx \right)^{k_2}\int_{\Omega} |u_n|^{r(x)-2}u_n (u_n-u) dx\nonumber.\\
\end{eqnarray}
Furthermore, utilizing H\"{o}lder's inequality and \eqref{cvg}, we can estimate:
\begin{eqnarray}\label{er7}
{\left|\int_{\Omega}| u_n|^{q(x)-2}u_n (u_n-u)dx\right|}&\leq& \int_{\Omega}| u_n|^{q(x)-1}| u_n-u| dx\notag\\
&\leq& C {\Big|{|  u_n|} ^{q(x)-1}\Big|} _{\frac{q(x)}{q(x)-1}}| u_n-u| _{q(x)}\notag\\&\leq& C\max\left\{\|u_n\|^{q^+-1}_{X}, \|u_n\|^{q^--1}_{X}\right\} | u_n-u| _{q(x)}.
\end{eqnarray}
Therefore, thanks to the convergence result \eqref{cvg}, we can deduce that
\begin{equation}\label{er8}
| u_n-u| _{q(x)}\to 0 \mbox{ as } n\to \infty.
\end{equation}
By combining the boundedness of $\{u_n\}$ in $X$ with the estimates \eqref{er7} and \eqref{er8}, we can conclude that
\begin{equation*}
\lim_{n\to \infty}\int_{\Omega}| u_n| ^{q(x)-2}u_n (u_n-u)dx=0.
\end{equation*}
As $\{u_n\}$ is bounded in $X$, there exist positive constants $c_1$ and $c_2$ such that
\begin{eqnarray}\label{c1c2}
c_1\leq \int_{\Omega}\frac{1}{q(x)} |u_n|^{q(x)}dx\leq c_2.
\end{eqnarray}
So, we have
\begin{equation}\label{cvgu0}
\left(\int_{\Omega}\frac{1}{q(x)} |u_n|^{q(x)}dx \right)^{k_1}\int_{\Omega} |u_n|^{q(x)-2}u_n (u_n-u) dx\to 0.
\end{equation}
Similarly, we obtain
\begin{equation}\label{cvgu1}
\lim_{n\to \infty}\left(\int_{\Omega}\frac{1}{r(x)} |u_n|^{r(x)}dx \right)^{k_2}\int_{\Omega} |u_n|^{r(x)-2}u_n (u_n-u) dx=0.
\end{equation}
By \eqref{er1}, we have
\[
\langle \mathcal{I}_{\lambda, \beta}'(u),u_n-u\rangle \to 0.
\]
Which means, based on equations \eqref{cvgu0} and \eqref{cvgu1},  that
\begin{equation}
\label{gj}
\left[a-b\left(\sigma_{p(x,y)}(u_n)\right)^{\gamma}\right]\int_{\Omega \times \Omega}\frac{|u_n(x)-u_n(y)|^{p(x,y)-2}(u_n(x)-u_n(y))((u_n(x)-u(x))-(u_n(y)-u(y)))}{|x-y|^{N+p(x,y)s(x,y)}}dxdy
\to 0.
\end{equation}
Since $\{u_n\}$ is bounded in $X$, passing to a subsequence, if necessary, we may assume that when $n\to \infty$
$$
\sigma_{p(x,y)}(u_n)=\int_{\Omega\times \Omega}\frac{1}{p(x,y)}\frac{|u_n(x)-u_n(y)|^{p(x,y)}}{|x-y|^{N+s(x,y){p(x,y)}}}dxdy\to t_0\geq 0.
$$
Considering two cases: $t_0 = 0$ and $t_0 > 0$. Now, let's proceed with a case analysis. First, if $t_0 = 0$, then the sequence $\{u_n\}$ converges strongly to $u = 0$ in $X$, and the proof is concluded. However, if $t_0 > 0$, we will further examine the two sub-cases below:
\\
{\bf Subcase 1.} If $t_0\neq \left(\frac ab\right)^\frac{1}{\gamma}$ then $a-b\left(\int_{\Omega\times \Omega}\frac{1}{p(x,y)}\frac{|u_n(x)-u_n(y)|^{p(x,y)}}{|x-y|^{N+s(x,y){p(x,y)}}}dxdy
\right)^{\gamma}\to 0$ is false, and there is no subsequence of \\$\left\{a-b\left(\int_{\Omega\times \Omega}\frac{1}{p(x,y)}\frac{|u_n(x)-u_n(y)|^{p(x,y)}}{|x-y|^{N+s(x,y){p(x,y)}}}dxdy\right)^{\gamma}\to 0\right\}$  converges to zero.  Thus, we can find a positive value $\delta > 0$
such that $$\left| a-b\left(\int_{\Omega\times \Omega}\frac{1}{p(x,y)}\frac{|u_n(x)-u_n(y)|^{p(x,y)}}{|x-y|^{N+s(x,y){p(x,y)}}}dxdy\right)^{\gamma}\right|  >\delta>0,$$
for sufficiently large $n$. As a result, we can conclude that the set
\begin{equation}
\left\{a-b\left(\int_{\Omega\times \Omega}\frac{1}{p(x,y)}\frac{|u_n(x)-u_n(y)|^{p(x,y)}}{|x-y|^{N+s(x,y){p(x,y)}}}dxdy\right)^{\gamma}\to 0\right\}  \mbox{ is
bounded}.
\end{equation}
{\bf Subcase 2.} If $t_0 =\left(\frac ab\right)^\frac{1}{\gamma}$, then $$a-b\left(\int_{\Omega\times \Omega}\frac{1}{p(x,y)}\frac{|u_n(x)-u_n(y)|^{p(x,y)}}{|x-y|^{N+s(x,y){p(x,y)}}}dxdy\right)^{\gamma}\to 0.$$
We define
$$
\varphi(u)=\frac{\lambda}{k_1+1}\left(\int_\Omega\frac{1}{q(x)} |u|^{q(x)}dx \right)^{k_1+1}+\frac{\beta}{k_2+1}\left(\int_\Omega\frac{1}{r(x)} |u|^{r(x)}dx \right)^{k_2+1},\; \mbox{for all} \; u\in X.
$$
Then
$$
\langle \varphi'(u),v\rangle=\lambda\left(\int_{\Omega}\frac{1}{q(x)} |u|^{q(x)}dx \right)^{k_1}\int_{\Omega} |u|^{q(x)-2}u v dx + \beta\left(\int_{\Omega}\frac{1}{r(x)} |u|^{r(x)}dx \right)^{k_2}\int_{\Omega} |u|^{r(x)-2}u v dx,\; \mbox{for all} \; v\in X.
$$
It follows that
\begin{eqnarray*}
\langle \varphi'(u_n)-\varphi'(u),v\rangle&=&\lambda\left[\left(\int_{\Omega}\frac{1}{q(x)} |u_n|^{q(x)}dx \right)^{k_1}\int_{\Omega} |u_n|^{q(x)-2}u_n v dx-\left(\int_{\Omega}\frac{1}{q(x)} |u|^{q(x)}dx \right)^{k_1}\int_{\Omega} |u|^{q(x)-2}u v dx\right]\\&\quad+&\beta\left[\left(\int_{\Omega}\frac{1}{r(x)} |u_n|^{r(x)}dx \right)^{k_2}\int_{\Omega} |u_n|^{r(x)-2}u_n v dx-\left(\int_{\Omega}\frac{1}{r(x)} |u|^{r(x)}dx \right)^{k_2}\int_{\Omega} |u|^{r(x)-2}u v dx\right].
\end{eqnarray*}
To complete our proof we require the following lemma.
\begin{lem}
\label{jkjk}
Suppose we have sequences $u_n$ and $u$ belonging to $X$ such that \eqref{cvg} is satisfied. Then, passing to a subsequence, if necessary, the following properties hold:
\begin{itemize}
  \item[(i)] $\lim_{n\to \infty}\left[\left(\int_{\Omega}\frac{1}{q(x)} |u_n|^{q(x)}dx \right)^{k_1}\int_{\Omega} |u_n|^{q(x)-2}u_n v dx-\left(\int_{\Omega}\frac{1}{q(x)} |u|^{q(x)}dx \right)^{k_1}\int_{\Omega} |u|^{q(x)-2}u v dx\right]=0$;
  \item[(ii)] $\lim_{n\to \infty}\left[\left(\int_{\Omega}\frac{1}{r(x)} |u_n|^{r(x)}dx \right)^{k_2}\int_{\Omega} |u_n|^{r(x)-2}u_n v dx-\left(\int_{\Omega}\frac{1}{r(x)} |u|^{r(x)}dx \right)^{k_2}\int_{\Omega} |u|^{r(x)-2}u v dx\right]=0$;
  \item[(iii)] $\langle \varphi'(u_n)-\varphi'(u),v\rangle\to 0,\;v\in X$.
\end{itemize}
\end{lem}
{\bf Proof.}  By \eqref{cvg}, we have  $u_n\to u$ in $L^{p(x)}(\O)$ which implies that
\begin{equation}
\label{cvgi}
|u_n|^{p(x)-2}u_n\to |u|^{p(x)-2}u \mbox{ in } L^{\frac{p(x)}{p(x)-1}}(\O).
\end{equation}
From \eqref{c1c2} we deduce that
\begin{eqnarray}
c_1^{k_1}\leq \left(\int_{\Omega}\frac{1}{q(x)} |u_n|^{q(x)}dx \right)^{k_1},\left(\int_{\Omega}\frac{1}{q(x)} |u|^{q(x)}dx \right)^{k_1}\leq c_2^{k_1}.
\end{eqnarray}
Due to H\"older's inequality, we have
\begin{eqnarray*}
&&\left|\left(\int_{\Omega}\frac{1}{q(x)} |u_n|^{q(x)}dx \right)^{k_1}\int_{\Omega} |u_n|^{q(x)-2}u_n v dx-\left(\int_{\Omega}\frac{1}{q(x)} |u|^{q(x)}dx \right)^{k_1}\int_{\Omega} |u|^{q(x)-2}u v dx\right|\\&=&\left|\left(\int_{\Omega}\frac{1}{q(x)} |u_n|^{q(x)}dx \right)^{k_1}\int_{\Omega} |u_n|^{q(x)-2}u_n v dx+\left(\int_{\Omega}\frac{1}{q(x)} |u|^{q(x)}dx \right)^{k_1}\int_{\Omega} |u|^{q(x)-2}u (-v) dx\right|\\&\leq&c_2^{k_1}
\left|\int_{\Omega} (|u_n|^{q(x)-2}u_n-|u|^{q(x)-2}u)v dx\right|\\&\leq& c_2^{k_1}\Big||u_n|^{q(x)-2}u_n-|u|^{q(x)-2}u\Big|_{\frac{q(x)}{q(x)-1}}|v|_{q(x)}\nonumber\\
&\leq& Cc_2^{k_1} \Big||u_n|^{q(x)-2}u_n-|u|^{q(x)-2}u\Big|_{\frac{q(x)}{q(x)-1}}\|v\|_X\to 0.
\end{eqnarray*}
By making a minor adjustment to the aforementioned proof, we can also establish assertion $(ii)$, but we will omit the specific details. As a result, by combining parts $(i)$ and $(ii)$, we can conclude assertion $(iii)$.
\\
Consequently, $\|  \varphi'(u_n)-\varphi'(u)\| _{X^*}\to 0$ and $\varphi'(u_n)\to\varphi'(u)$. \qed

We are now able to conclude the proof of Subcase $2$. Utilizing Lemma \ref{jkjk} and taking into account the fact that $\langle \mathcal{I}'_{\lambda,\beta}(u),v\rangle=\left[a-b\left(\sigma_{p(x,y)}(u)\right)^{\gamma}\right]\int_{\Omega\times \Omega}\frac{|u(x)-u(y)|^{p(x,y)-2}(u(x)-u(y))(v(x)-v(y))}{|x-y|^{N+p(x,y)s(x,y)}}dxdy  -\langle \varphi'(u),v\rangle$,\\$\langle \mathcal{I}'_{\lambda,\beta}((u_n),v\rangle\to 0$
  and $a-b\left(\int_{\Omega\times \Omega}\frac{1}{p(x,y)}\frac{|u_n(x)-u_n(y)|^{p(x,y)}}{|x-y|^{N+s(x,y){p(x,y)}}}dxdy\right)^{\gamma}\to 0$, then we can infer that $\varphi'(u_n)\to 0\;(n\to\infty)$, i.e.,~$$
\langle \varphi'(u),v\rangle=\lambda\left(\int_{\Omega}\frac{1}{q(x)} |u|^{q(x)}dx \right)^{k_1}\int_{\Omega} |u|^{q(x)-2}u v dx + \beta\left(\int_{\Omega}\frac{1}{r(x)} |u|^{r(x)}dx \right)^{k_2}\int_{\Omega} |u|^{r(x)-2}u v dx,\; \mbox{for all} \; v\in X,
$$
 and therefore
$$
\lambda\left(\int_{\Omega}\frac{1}{q(x)} |u|^{q(x)}dx \right)^{k_1}|u(x)|^{q(x)-2}u(x)  + \beta \left(\int_{\Omega}\frac{1}{r(x)} |u|^{r(x)}dx \right)^{k_2}|u(x)|^{r(x)-2}u(x) =0\mbox{ for a.e.} \; x \in \Omega.
$$
By invoking the fundamental lemma of the variational method (see~\cite{Willem}), we can conclude that $u = 0$. Hence,
\begin{gather*}
\varphi(u_n)=\frac{\lambda}{k_1+1}\left(\int_\Omega\frac{1}{q(x)} |u_n|^{q(x)}dx \right)^{k_1+1}+\frac{\beta}{k_2+1}\left(\int_\Omega\frac{1}{r(x)} |u_n|^{r(x)}dx \right)^{k_2+1}\\\to \frac{\lambda}{k_1+1}\left(\int_\Omega\frac{1}{q(x)} |u|^{q(x)}dx \right)^{k_1+1}+\frac{\beta}{k_2+1}\left(\int_\Omega\frac{1}{r(x)} |u|^{r(x)}dx \right)^{k_2+1}=0.
\end{gather*}
Hence, we can deduce that
\begin{eqnarray*}
\mathcal{I}_{\lambda, \beta}(u_n)&=&a\int_{\Omega\times \Omega}\frac{1}{p(x,y)}\frac{|u_n(x)-u_n(y)|^{p(x,y)}}{|x-y|^{N+s(x,y){p(x,y)}}} \,dx\,dy-\frac{b}{\gamma+1}\left(\int_{\Omega\times \Omega}\frac{1}{p(x,y)}\frac{|u_n(x)-u_n(y)|^{p(x,y)}}{|x-y|^{N+s(x,y){p(x,y)}}} \,dx\,dy\right)^{\gamma+1}\\& -&\frac{\lambda}{k_1+1}\left(\int_\Omega\frac{1}{q(x)} |u|^{q(x)}dx \right)^{k_1+1}-\frac{\beta}{k_2+1}\left(\int_\Omega\frac{1}{r(x)} |u|^{r(x)}dx \right)^{k_2+1}\\&&\to \frac{\gamma a^\frac{\gamma+1}{\gamma}}{(\gamma+1)b^\frac{1}{\gamma}}.
\end{eqnarray*}
Therefore, we have reached a contradiction since $\mathcal{I}_{\lambda, \beta}(u_n)\to c\in\left(0, \frac{a^\frac{\gamma+1}{\gamma}}{ b^\frac{1}{\gamma}+\frac{b^\frac{1}{\gamma}}{\gamma}}\right)$. \\Then $a-b\left(\int_{\Omega\times \Omega}\frac{1}{p(x,y)}\frac{|u_n(x)-u_n(y)|^{p(x,y)}}{|x-y|^{N+s(x,y){p(x,y)}}} \,dx\,dy\right)^{\gamma}\to 0$ is not true. Similarly to Subcase 1, we can argue as follows:
\begin{equation*}
\left\{a-b\left(\int_{\Omega\times \Omega}\frac{1}{p(x,y)}\frac{|u_n(x)-u_n(y)|^{p(x,y)}}{|x-y|^{N+s(x,y){p(x,y)}}} \,dx\,dy\right)^{\gamma}\to 0\right\}  \mbox{ is
bounded}.
\end{equation*}
So, combining the two cases discussed above, we can conclude that:
$$
\int_{\Omega \times \Omega}\frac{|u_n(x)-u_n(y)|^{p(x,y)-2}(u_n(x)-u_n(y))((u_n(x)-u(x))-(u_n(y)-u(y)))}{|x-y|^{N+p(x,y)s(x,y)}}dxdy\to 0.
$$
Therefore, by invoking the $(S_+)$ condition and Proposition \ref{lem 3.2}, we conclude that $\| u_n\|_X \to\| u\|_X $ as $n\to \infty$, which implies that $\mathcal{I}_{\lambda, \beta}$ satisfies the $(PS)_c$ condition. Hence, the proof is now complete.\qed
\subsection{\bf{The $(PS)_c$ condition for $\mathbf{a=0}$}}
\begin{lem}\label{lem2}
Assuming that \eqref{cond(1.4)a=0} is valid, then the functional $\mathcal{I}_{\lambda, \beta}$ satisfies the Palais-Smale condition at all level $c\in \R$.
\end{lem}
{\bf Proof.} Let $\{u_n\}$ be a $(PS)_c$ sequence of $\mathcal{I}_{\lambda, \beta}$, that is
\begin{eqnarray}\label{err1}
\mathcal{I}_{\lambda, \beta}(u_n) \to c, \quad \mathcal{I}'_{\lambda, \beta}(u_n) \to 0 \mbox{ in } X^*, \quad n \to \infty,
\end{eqnarray}
where $X^*$ is the dual space of $X$.
\\
{\bf Step 1.} We will prove that $\{u_n\}$ is bounded in $X$. Let us assume by contradiction that $\{u_n\}$ is unbounded in $X$. Without loss of generality, we can assume that $\|u_n\|_X > 1$ for all $n$. Take $$\theta<\min\left\{\frac{(\gamma+1)(p^-)^{\gamma+1}}{(p^+)^{\gamma}},\frac{(k_1+1)(q^-)^{k_1+1}}{(q^+)^{k_1}},\frac{(k_2+1)(r^-)^{k_2+1}}{(r^+)^{k_2}}\right\},$$
then, we have
\begin{eqnarray}\label{err2.}
& &c+1+\|u_n\|_X  \geq  \mathcal{I}_{\lambda, \beta}(u_n) - \frac{1}{\theta}\langle \mathcal{I}_{\lambda, \beta}'(u_n),u_n\rangle \notag\\
& = & \left[-\frac{b}{\gamma+1}\left(\sigma_{p(x,y)}(u_n)\right)^{\gamma+1} -\frac{\lambda}{k_1+1}\left(\int_\Omega\frac{1}{q(x)} |u_n|^{q(x)}dx \right)^{k_1+1}-\frac{\beta}{k_2+1}\left(\int_\Omega\frac{1}{r(x)} |u_n|^{r(x)}dx \right)^{k_2+1}\right]\notag\\
&\qquad -\frac{1}{\theta}&\Bigg(\left[-b\left(\sigma_{p(x,y)}(u_n)\right)^{\gamma}\right]\int_{\Omega\times \Omega}\frac{|u_n(x)-u_n(y)|^{p(x,y)}}{|x-y|^{N+p(x,y)s(x,y)}}dxdy -\lambda\left(\int_{\Omega}\frac{1}{q(x)} |u_n|^{q(x)}dx \right)^{k_1}\int_{\Omega} |u_n|^{q(x)}dx\notag\\ & \qquad-& \beta\left(\int_{\Omega}\frac{1}{r(x)} |u_n|^{r(x)}dx \right)^{k_2}\int_{\Omega} |u_n|^{r(x)}dx\Bigg)\notag\\
&\geq&b\left(\frac{1}{\theta( p_+)^\gamma}-\frac{1}{(\gamma+1)(p^-)^{\gamma+1}}\right)\left(\int_{\Omega\times \Omega}\frac{|u_n(x)-u_n(y)|^{p(x,y)}}{|x-y|^{N+p(x,y)s(x,y)}}dxdy\right)^{\gamma+1}\notag\\
&\qquad+&\lambda\left(\frac{1}{\theta(q^+)^{k_{1}}}-\frac{1}{(k_1+1)(q^-)^{k_1+1}}\right)\left(\int_{\Omega} |u_n|^{q(x)}dx \right)^{k_1+1}+\beta\left(\frac{1}{\theta(r^+)^{k_{2}}}-\frac{1}{(k_2+1)(r^-)^{k_2+1}}\right)\left(\int_{\Omega} |u_n|^{r(x)}dx \right)^{k_2+1}
\notag\\
&\geq&b\left(\frac{1}{\theta( p_+)^\gamma}-\frac{1}{(\gamma+1)(p^-)^{\gamma+1}}\right)\|u_n\|^{(\gamma+1)p^-}_X
+\lambda\left(\frac{1}{\theta(q^+)^{k_{1}}}-\frac{1}{(k_1+1)(q^-)^{k_1+1}}\right)\left(\int_{\Omega} |u_n|^{q(x)}dx \right)^{k_1+1}\notag\\&\qquad +&\beta\left(\frac{1}{\theta(r^+)^{k_{2}}}-\frac{1}{(k_2+1)(r^-)^{k_2+1}}\right)\left(\int_{\Omega} |u_n|^{r(x)}dx \right)^{k_2+1}.
\end{eqnarray}
For simplicity, let's denote
\begin{eqnarray}\label{err3a=0}
A_1 &=& b\left(\frac{1}{\theta( p^+)^\gamma}-\frac{1}{(\gamma+1)(p^-)^{\gamma+1}}\right),\nonumber\\
A_2 &=& \lambda\left(\frac{1}{\theta(q^+)^{k_{1}}}-\frac{1}{(k_1+1)(q^-)^{k_1+1}}\right),\nonumber\\
A_3 &=& \beta\left(\frac{1}{\theta(r^+)^{k_{2}}}-\frac{1}{(k_2+1)(r^-)^{k_2+1}}\right).
\end{eqnarray}
Using \eqref{err2.} and \eqref{err3a=0}, we can write
\begin{eqnarray}\label{er3a=0}
A_1\|u_n\|^{(\gamma+1)p^-}_X&\leq&\left\{
 \begin{array}{ll}
c+1+\|u_n\|_X-A_3\left(\int_{\Omega} |u_n|^{r(x)}dx \right)^{k_2+1} ,\quad \mbox{ if } \lambda>0, \beta<0.\nonumber\\\\

c+1+\|u_n\|_X-A_2\left(\int_{\Omega} |u_n|^{q(x)}dx \right)^{k_1+1},\quad \mbox{ if } \lambda<0, \beta>0.\nonumber\\\\

c+1+\|u_n\|_X,\quad \mbox{ if } \lambda<0, \beta<0.
\end{array}
\right.
\\
&\leq&\left\{
 \begin{array}{ll}
c+1+\|u_n\|_X-A_3\|u_n\|_X^{r^+(k_2+1)} ,\quad \mbox{ if } \lambda>0, \beta<0. \\\\

c+1+\|u_n\|_X-A_2\|u_n\|_X^{q^+(k_1+1)},\quad \mbox{ if } \lambda<0, \beta>0. \\\\

c+1+\|u_n\|_X,\quad \mbox{ if } \lambda<0, \beta<0.
\end{array}
\right.
\end{eqnarray}
It follows from \eqref{cond(1.4)a=0} and \eqref{er3a=0} that $\{u_n\}$ is bounded in $X$.
\\
{\bf Step 2.}  We will now demonstrate that the sequence $\{u_n\}$ possesses a convergent subsequence in the space $X$. According to Theorem \ref{prp 3.3}, the embedding $X \hookrightarrow L^{\tau(x)}(\Omega)$ is compact where $1\leq \tau(x)<p^*_{s}(x)$. Since $X$ is a reflexive Banach space, passing, if necessary, to a subsequence, there exists $u\in X$ such that
\begin{equation}
\label{cvgg}
u_n\rightharpoonup u\mbox{ in } X,\; u_n \to u \mbox{ in } L^{\tau(x)}(\Omega),\;\; u_n(x)\to u(x), \mbox{ a.e. in } \Omega.
\end{equation}
From \eqref{deriv-fun}, we find that
\begin{eqnarray} \label{er6a=0}
& &\langle \mathcal{I}_{\lambda, \beta}'(u),u_n-u\rangle\nonumber\\ &=& -b\left(\sigma_{p(x,y)}(u_n)\right)^{\gamma}\int_{\Omega\times \Omega}\frac{|u_n(x)-u_n(y)|^{p(x,y)-2}(u_n(x)-u_n(y))((u_n(x)-u(x))-(u_n(y)-u(y)))}{|x-y|^{N+p(x,y)s(x,y)}}dxdy \nonumber\\&\quad-&\lambda\left(\int_{\Omega}\frac{1}{q(x)} |u_n|^{q(x)}dx \right)^{k_1}\int_{\Omega} |u_n|^{q(x)-2}u_n (u_n-u) dx -\beta \left(\int_{\Omega}\frac{1}{r(x)} |u_n|^{r(x)}dx \right)^{k_2}\int_{\Omega} |u_n|^{r(x)-2}u_n (u_n-u) dx\nonumber\\
\end{eqnarray}
So, we have
\begin{equation}\label{cvguu0}
\left(\int_{\Omega}\frac{1}{q(x)} |u_n|^{q(x)}dx \right)^{k_1}\int_{\Omega} |u_n|^{q(x)-2}u_n (u_n-u) dx\to 0.
\end{equation}
Similarly, we obtain
\begin{equation}\label{cvguu1}
\lim_{n\to \infty}\left(\int_{\Omega}\frac{1}{r(x)} |u_n|^{r(x)}dx \right)^{k_2}\int_{\Omega} |u_n|^{r(x)-2}u_n (u_n-u) dx=0.
\end{equation}
By \eqref{err1}, we have
\[
\langle \mathcal{I}_{\lambda, \beta}'(u),u_n-u\rangle \to 0.
\]
So, based on the expressions \eqref{cvguu0} and \eqref{cvguu1}, we can conclude that \eqref{er6a=0} leads to the following implications:
\begin{equation}
\label{gj}
-b\left(\sigma_{p(x,y)}(u_n)\right)^{\gamma}\int_{\Omega \times \Omega}\frac{|u_n(x)-u_n(y)|^{p(x,y)-2}(u_n(x)-u_n(y))((u_n(x)-u(x))-(u_n(y)-u(y)))}{|x-y|^{N+p(x,y)s(x,y)}}dxdy
\to 0.
\end{equation}
Since $\{u_n\}$ is bounded in $X$ and $b>0$, we have
\begin{equation*}
\left\{-b\left(\int_{\Omega\times \Omega}\frac{1}{p(x,y)}\frac{|u_n(x)-u_n(y)|^{p(x,y)}}{|x-y|^{N+s(x,y){p(x,y)}}} \,dx\,dy\right)^{\gamma}\to 0\right\}  \mbox{ is
bounded}.
\end{equation*}
Therefore, we can conclude from the two aforementioned cases that
$$
\int_{\Omega \times \Omega}\frac{|u_n(x)-u_n(y)|^{p(x,y)-2}(u_n(x)-u_n(y))((u_n(x)-u(x))-(u_n(y)-u(y)))}{|x-y|^{N+p(x,y)s(x,y)}}dxdy\to 0.
$$
Therefore, by utilizing the $(S_+)$ condition and Proposition \ref{lem 3.2}, we can deduce that $\| u_n\|_X \to\| u\|_X $ as $n\to \infty$, indicating that $\mathcal{I}_{\lambda,\beta}$ satisfies the $(PS)_c$ condition. This concludes the proof. \qed
\section{Proof of Theorem \ref{theo1.1}}
In this part, one is addressed in proving Theorem \ref{theo1.1} by applying the mountain pass theorem, see \cite{Willem}.
\begin{lem}\label{lemme1}
Assume that \eqref{cond-pqr} holds. Then there exist $\rho > 0$ and
$\alpha > 0$ such that $\mathcal{I}_{\lambda,\beta}(u) \geq \alpha > 0$, for any $u \in X$ with $\| u\|_X =\rho$.
\end{lem}
{\bf Proof.} Let $u\in X$ with $\|u\|_X < 1$. From \eqref{er1}, Lemma \ref{lem 3.1} and Sobolev immersions, we get
\begin{eqnarray*}\label{er10}
\mathcal{I}_{\lambda, \beta}(u) & = & a\sigma_{p(x,y)}(u)-\frac{b}{\gamma+1}\left(\sigma_{p(x,y)}(u)\right)^{\gamma+1} -\frac{\lambda}{k_1+1}\left(\int_\Omega\frac{1}{q(x)} |u|^{q(x)}dx \right)^{k_1+1}-\frac{\beta}{k_2+1}\left(\int_\Omega\frac{1}{r(x)} |u|^{r(x)}dx \right)^{k_2+1} \notag\\
& \geq & \frac{a}{p^+}\|u\|^{p^+}_X - \frac{b}{(p^-)^{\gamma+1}(\gamma+1)}\|u\|^{p^-(\gamma+1)}_X-\frac{\lambda}{k_1+1}\frac{1}{(q^-)^{k_1+1}}\left(\int_\Omega|u|^{q(x)}dx \right)^{k_1+1}-\frac{\beta}{k_2+1}\frac{1}{(r^-)^{k_2+1}}\left(\int_\Omega|u|^{r(x)}dx \right)^{k_2+1}\notag\\
& \geq & \frac{a}{p^+}\|u\|^{p^+}_X - \frac{b}{(p^-)^{\gamma+1}(\gamma+1)}\|u\|^{p^{-}(\gamma+1)}_X-\frac{\lambda}{k_1+1}\frac{C_1^{q^-(k_1+1)}}{(q^-)^{k_1+1}}\|u\|_X^{q^-(k_1+1)}
-\frac{\beta}{k_2+1}\frac{C_2^{r^-(k_2+1)}}{(r^-)^{k_2+1}}\|u\|_X^{r^-(k_2+1)}\\
&\geq&\|u\|^{p^+}_X\left(\frac{a}{p^+}- \frac{b}{(p^-)^{\gamma+1}(\gamma+1)}\|u\|^{p^-(\gamma+1)-p^{+}}_X-\frac{\lambda}{k_1+1}\frac{C_1^{q^-(k_1+1)}}{(q^-)^{k_1+1}}\|u\|_X^{q^-(k_1+1)-p^{+}}
-\frac{\beta}{k_2+1}\frac{C_2^{r^-(k_2+1)}}{(r^-)^{k_2+1}}\|u\|_X^{r^-(k_2+1)-p^{+}}\right).
\end{eqnarray*}
Hence, based to the fact that $\|u\|_{X}<1$ and  $p$ satisfies the condition\eqref{cond-pqr}, we infer the result.\qed

\begin{lem}\label{lemme2}
Assume that the conditions \ref{hypot:S},\ref{hypot:P},  \eqref{cond-pqr} hold. Then there exists $e\in X$ with $\| e\|_X >\rho$ (where $\rho$ is given by  Lemma \ref{lemme1}) such that $\mathcal{I}_{\lambda, \beta}(e) < 0$.
\end{lem}
{\bf Proof.} Let $\phi_0 \in C_0^\infty(\O)$. According to the condition \eqref{cond-pqr}, for $t>1$ large enough, we have
\begin{eqnarray*}\label{er4}
\mathcal{I}_{\lambda, \beta}(t\phi_0) & \leq & \frac{at^{p^+}}{p^-}\|\phi_0\|^{p^+}_X-\frac{bt^{p^-(\gamma+1)}}{(\gamma+1)(p^+)^{\gamma+1}}\|\phi_0\|^{p^-(\gamma+1)}_X -\frac{\lambda t^{q^-(k_1+1)}}{(k_1+1)(q^+)^{k_1+1}}\left(\int_\Omega |\phi_0|^{q(x)}dx \right)^{k_1+1}\notag\\&-&\frac{\beta t^{r^-(k_2+1)}}{(k_2+1)(r^+)^{k_2+1}}\left(\int_\Omega |\phi_0|^{r(x)}dx \right)^{k_2+1}.
\end{eqnarray*}
If the condition \eqref{cond-pqr} holds, then $\mathcal{I}_{\lambda, \beta}(t\phi_0) \to -\infty$ as $t\to \infty$. So, for some $t_0>1$ large enough, we deduce that $\|t_0\phi_0\|_X > \rho$ and $\mathcal{I}_{\lambda, \beta}(t_0\phi_0) < 0$. Choosing $e = t_0\phi_0$, the proof of Lemma \ref{lemme2} is completed.\qed
\subsection*{Proof of Theorem \ref{theo1.1}.}
 It follows from Lemmas \ref{lem2}, \ref{lemme1}, \ref{lemme2} and the fact that $\mathcal{I}_{\lambda, \beta}(0)=0$, $\mathcal{I}_{\lambda, \beta}$ satisfies all conditions of the mountain pass theorem \cite{Willem}. Thus, problem \eqref{main} admits a nontrivial weak solution.
\qed

\section{Proof of Theorem \ref{theo1.2}}
 Since $X$ is a reflexive and separable Banach space, so there exist $e_{i}\in X$ and $e_{i}^*\in X^*$ such that $\langle e_i, e^*_{j}\rangle =\delta_{ij}$ where $\delta_{}$ means the Kronecker symbol.\\
We denote $$X_{i}=\overline {span\{e_i, i=1,2,\cdots\}},\; X^*_{i}=\overline {span\{e^*_{i}, i=1,2,\cdots\}}.$$
Now, we consider $X_i=\{e_i\}$ and let denote
$$X=\overline{\oplus_{i=1}^{\infty} X_i}, \quad
Y_{k}=\oplus_{i=1}^{k} X_i, \quad
Z_{k}=\overline{\oplus_{i=k}^{\infty} X_i}.$$
\begin{theorem}[{Fountain Theorem, see~\cite{Willem}}] \label{fountain}
 Let $ X_{0} $ be a Banach space with the norm $ \|\cdot\|_{X_{0}} $ and  let $ X_i $
be a sequence of subspace of $ X_{0} $ with $ dim X_i < \infty $ for each
$ i\in N $. In addition, set
\begin{equation*}
 X_{0}=\overline{\oplus_{i=1}^{\infty} X_i}, \quad
Y_{k}=\oplus_{i=1}^{k} X_i, \quad
Z_{k}=\overline{\oplus_{i=k}^{\infty} X_i}
 \end{equation*}
For each even functional $ J \in C^{1}(X_{0},\mathbb{R})$ and for each  $ k \in \mathbb{N} $, we suppose that  there exists
$ \rho_{k} > \gamma_{k}> 0 $ such that
\begin{itemize}
\item[(1)] $a_{k}:= \max_{u\in Y_{k},\| u \|_{X_{0}} = \rho_{k}} J(u)\leq  0$,

\item[(2)] $b_{k}:=\inf_{u\in Z_{k}, \| u \|_{X_{0}} =\gamma_{k}} J(u) \to +\infty$, $k \to  +\infty$,

\item[(3)] The functional $J$ satisfies the $(\text{PS})_{\text{c}}$ condition for every $c > 0$.
\end{itemize}
 Then $ J $ admits an unbounded sequence of critical values.
\end{theorem}

To prove our result, we will use the Fountain theorem \ref{fountain}. So, this proof is divided  in several lemmas given as follows
\begin{lem}(see \cite{BT(CVEE)})
If $q(x), r(x)\in C_{+}(\overline{\Omega})$ satisfying $1\leq q(x), r(x)<p_{s}^{*}(x), \forall x\in \overline{\Omega}$ and let denote by
$$\xi_{k}=\sup\{|u|_{q(x)}, \|u\|_{X}=1,\; u\in Z_{k}\},\;\;\xi'_{k}=\sup\{|u|_{r(x)}, \|u\|_{X}=1,\; u\in Z_{k}\}.$$
Then  $$\lim_{k\to \infty}\xi_{k}= 0,\; \lim_{k\to \infty}\xi'_{k}= 0.$$
\end{lem}
\begin{lem}
The functional $I_{\lambda, \beta}$ verifies the following property
$a_{k}:= \max_{u\in Y_{k},\| u \| = \rho_{k}} I_{\lambda, \beta}(u)\leq  0$,
where the space $Y_{k}$ is given in the Theorem \ref{fountain}.
\end{lem}{\bf Proof.}
Let $\lambda ,\beta >0.$ Since $Y_{k}=\oplus_{i=1}^{k} X_i$, then $dim Y_k<\infty$ or all norms are equivalent in the finite dimensional space. Let now $u\in Y_k$ such that $\|u\|_{X}>1$,  thus we have
\begin{eqnarray*}
\mathcal{I}_{\lambda, \beta}(u) & = & a\sigma_{p(x,y)}(u)-\frac{b}{\gamma+1}\left(\sigma_{p(x,y)}(u)\right)^{\gamma+1} -\frac{\lambda}{k_1+1}\left(\int_\Omega\frac{1}{q(x)} |u|^{q(x)}dx \right)^{k_1+1}-\frac{\beta}{k_2+1}\left(\int_\Omega\frac{1}{r(x)} |u|^{r(x)}dx \right)^{k_2+1} \notag\\
& \leq & \frac{a}{p^-}\|u\|^{p^+}_X - \frac{b}{(p^+)^{\gamma+1}(\gamma+1)}\|u\|^{p^-(\gamma+1)}_X-\frac{\lambda}{k_1+1}\left(\int_\Omega\frac{1}{q(x)} |u|^{q(x)}dx \right)^{k_1+1}-\frac{\beta}{k_2+1}\left(\int_\Omega\frac{1}{r(x)} |u|^{r(x)}dx \right)^{k_2+1}\\
& \leq & \frac{a}{p^-}\|u\|^{p^+}_X - \frac{b}{(p^+)^{\gamma+1}(\gamma+1)}\|u\|^{p^-(\gamma+1)}_X-\frac{\lambda}{(k_1+1){(q^{+})}^{k_{1}+1}}\left(\int_\Omega|u|^{q(x)}dx \right)^{k_1+1}-\frac{\beta}{(k_2+1){(r^{+})}^{k_{2}+1}}\left(\int_\Omega|u|^{r(x)}dx \right)^{k_2+1}
\end{eqnarray*}
Hence, By using the embeddings $L^{q(x)}   \hookrightarrow X$ and $L^{r(x)}   \hookrightarrow X$ (see theorem \ref{prp 3.3}) and based to the inequality \eqref{cond-pqr}, we infer that $a_{k}:= \max_{u\in Y_{k},\| u \| = \rho_{k}} I_{\lambda, \beta}(u)\leq  0$ \qed
\begin{lem}
The functional $I_{\lambda, \beta}$ verifies the following property
$b_{k}:=\inf_{u\in Z_{k}, \| u \| =\gamma_{k}} J(u) \to +\infty$, $k \to  +\infty$,
where the space $Z_{k}$ is given in the Theorem \ref{fountain}.
\end{lem}
{\bf Proof.}
Let $u\in Z_k$ with $\|u\|_{X}<1$. So, we have
\begin{eqnarray*}
\mathcal{I}_{\lambda, \beta}(u) & = & a\sigma_{p(x,y)}(u)-\frac{b}{\gamma+1}\left(\sigma_{p(x,y)}(u)\right)^{\gamma+1} -\frac{\lambda}{k_1+1}\left(\int_\Omega\frac{1}{q(x)} |u|^{q(x)}dx \right)^{k_1+1}-\frac{\beta}{k_2+1}\left(\int_\Omega\frac{1}{r(x)} |u|^{r(x)}dx \right)^{k_2+1} \notag\\
& \geq & \frac{a}{p^+}\|u\|^{p^+}_X - \frac{b}{(p^-)^{\gamma+1}(\gamma+1)}\|u\|^{p^-(\gamma+1)}_X-\frac{\lambda}{k_1+1}\frac{1}{(q^-)^{k_1+1}}\left(\int_\Omega|u|^{q(x)}dx \right)^{k_1+1}-\frac{\beta}{k_2+1}\frac{1}{(r^-)^{k_2+1}}\left(\int_\Omega|u|^{r(x)}dx \right)^{k_2+1}\notag\\
\end{eqnarray*}
So, we obtain
 $$\mathcal{I}_{\lambda, \beta}(u)\geq \begin{cases}
\frac{a}{p^+}\|u\|^{p^+}_{X} - \frac{b}{(p^-)^{\gamma+1}(\gamma+1)}\|u\|^{p^{-}(\gamma+1)}_{X}-\lambda C_{q}-\beta C_{r},\;\text{if}\; |u|_{q(x)}<1,\; \text{and}\; |u|_{r(x)}<1,\\
\frac{a}{p^+}\|u\|^{p^+}_{X} - \frac{b}{(p^-)^{\gamma+1}(\gamma+1)}\|u\|^{p^-(\gamma+1)}_X-\lambda C'_{q}[\xi_{k}\|u\|_X]^{q^-(k_{1}+1)}
-\beta C_{r},\; \text{if} |u|_{q(x)}>1\;\text{and}\;  |u|_{r(x)}<1,\\
\frac{a}{p^+}\|u\|^{p^+}_{X} - \frac{b}{(p^-)^{\gamma+1}(\gamma+1)}\|u\|^{p^-(\gamma+1)}_{X}-\lambda C_{q}-\beta C'_{r}[\xi'_{k}\|u\|_X]^{r^-(k_{2}+1)},\;
\text{if} |u|_{q(x)}<1\;\text{and}\;  |u|_{r(x)}>1\\
\frac{a}{p^+}\|u\|^{p^+}_X - \frac{b}{(p^-)^{\gamma+1}(\gamma+1)}\|u\|^{p^-(\gamma+1)}_{X}-\lambda C'_{q}[\xi_{k}\|u\|_X]^{q^-(k_{1}+1)}-\beta C'_{r}[\xi'_{k}\|u\|_X]^{r^-(k_{2}+1)},\;
\text{if}\; |u|_{q(x)}>1\;\text{and}\;  |u|_{r(x)}>1,
\end{cases}
$$
where, $$C_q=
\frac{1}{k_1+1}\frac{1}{(q^-)^{k_1+1}}C^{q^{-}(k_{1}+1)},\;  C'_q=
\frac{1}{k_1+1}\frac{1}{(q^-)^{k_1+1}}C^{q^{+}(k_{1}+1)}\;$$ and $$C_r=\frac{1}{k_2+1}\frac{1}{(r^-)^{k_2+1}}C^{r^{-}(k_{2}+1)},\; C'_r=\frac{1}{k_2+1}\frac{1}{(r^-)^{k_2+1}}C^{r^{+}(k_{2}+1)}.$$
Hence, we have $$\mathcal{I}_{\lambda, \beta}(u)
\geq\frac{a}{p^-}\|u\|^{p^+}_X - \frac{b}{(p^+)^{\gamma+1}(\gamma+1)}\|u\|^{p^-(\gamma+1)}_X-\lambda C''_{q}[\xi_{k}\|u\|_X]^{q^-(k_1+1)}-\beta C''_{r}[\xi'_{k}\|u\|_X]^{r^-(k_2+1)}-C_{\lambda, \beta}.$$
So, based to the fact that $\lim_{k\to \infty}\xi_{k}=0$ and $\lim_{k\to \infty}\xi'_{k}=0$, we can deduce that for $k$ sufficiently large, we have $\xi_{k}<1$ and $\xi'_{k}<1$. Thus, we have
$b_{k}:=\inf_{u\in Z_{k}, \| u \| =\gamma_{k}} J(u) \to +\infty$, $k \to  +\infty$  since we have taken $\|u\|_{X}<1.$
\qed

\textbf{Proof of Theorem \ref{theo1.2}}
We have that $\mathcal{I}_{\lambda, \beta}(u)$ belongs to $C^{1}(X,\mathbb{R})$, even functional and verifies the Palais-smale condition. Moreover, we have
$$a_{k}:= \max_{u\in Y_{k},\| u \| = \rho_{k}} I_{\lambda, \beta}(u)\leq  0$$ and
$$ b_{k}:=\inf_{u\in Z_{k}, \| u \| =\gamma_{k}} J(u) \to +\infty,$$ for $k \to  +\infty.$
Then, by using the Fountain theorem, we deduce that $\mathcal{I}_{\lambda, \beta}(u)$ admits an unbounded sequence of critical points.
\subsection{\bf Proof of Theorem \ref{symm.moun.pass}.}
To prove Theorem \ref{symm.moun.pass}, we shall use the following symmetric mountain pass theorem in \cite{HZCR}:
\begin{theorem} \label{th45}(\cite{HZCR}).
 Let $E$  be a real infinite dimensional Banach space and $I\in C^1(E)$ satisfying the Palais-Smale condition. Suppose $E=E^-\oplus E^+$, where
 $E^-$ is finite dimensional, and assume the following conditions:
\begin{enumerate}
  \item $I$ is even and $I(0)=0$;
 \item there exist $\alpha >0$ and $\rho>0$ such that $I(u)\geq \alpha$ for any $u\in E^+$ with $\|u\|=\rho$;
\item for any finite dimensional subspace $W\subset E$ there is $R= R(W)$ such that $I(u)\leq 0$ for $u\in W$, $\|u\|\geq R$;
\end{enumerate}
then, $I$ possesses an unbounded sequence of critical values.
\end{theorem}

\begin{lem}\label{lemme1K}
Assume that \eqref{cond(1.4)a=0} holds. Then there exist $\rho > 0$ and
$\alpha > 0$ such that $\mathcal{I}_{\lambda,\beta}(u) \geq \alpha > 0$, for any $u \in X$ with $\| u\|_X =\rho$.
\end{lem}
{\bf Proof.} Let $u\in X$ with $\|u\|_X =\rho\in(0,1)$. Following \cite{ABMS}, let denote by $\lambda_{1}(x)$ the eigenvalue related to our operator. By using the Sobolev immersions, we get
\begin{eqnarray*}
\mathcal{I}_{\lambda, \beta}(u) & = & -\frac{b}{\gamma+1}\left(\sigma_{p(x,y)}(u)\right)^{\gamma+1} -\frac{\lambda}{k_1+1}\left(\int_\Omega\frac{1}{q(x)} |u|^{q(x)}dx \right)^{k_1+1}-\frac{\beta}{k_2+1}\left(\int_\Omega\frac{1}{r(x)} |u|^{r(x)}dx \right)^{k_2+1} \notag
\\
&\geq &\left\{
 \begin{array}{ll}
-\frac{\beta}{k_2+1}\left(\int_\Omega\frac{1}{r(x)} |u|^{r(x)}dx \right)^{k_2+1}- \frac{b}{(p^-)^{\gamma+1}(\gamma+1)}\|u\|^{p^-(\gamma+1)}_X-\frac{\lambda}{k_1+1}\left(\int_\Omega\frac{1}{q(x)} |u|^{q(x)}dx \right)^{k_1+1},\\\quad \mbox{ if } \lambda>0, \beta<0. \\\\

-\frac{\lambda}{k_1+1}\left(\int_\Omega\frac{1}{q(x)} |u|^{q(x)}dx \right)^{k_1+1}- \frac{b}{(p^-)^{\gamma+1}(\gamma+1)}\|u\|^{p^-(\gamma+1)}_X -\frac{\beta}{k_2+1}\left(\int_\Omega\frac{1}{r(x)} |u|^{r(x)}dx \right)^{k_2+1},\\\quad \mbox{ if } \lambda<0, \beta>0. \\\\

-\frac{\lambda}{k_1+1}\left(\int_\Omega\frac{1}{q(x)} |u|^{q(x)}dx \right)^{k_1+1}-\frac{\beta}{k_2+1}\left(\int_\Omega\frac{1}{r(x)} |u|^{r(x)}dx \right)^{k_2+1}- \frac{b}{(p^-)^{\gamma+1}(\gamma+1)}\|u\|^{p^-(\gamma+1)}_X,\\\quad \mbox{ if } \lambda<0, \beta<0.
\end{array}
\right.
\\
&\geq &\left\{
 \begin{array}{ll}
-\frac{\beta\lambda_1(x)}{(r^+)^{k_2+1}(k_2+1)}
\|u\|^{r^+(k_2+1)}- \frac{b}{(p^-)^{\gamma+1}(\gamma+1)}\|u\|^{p^-(\gamma+1)}_X
-\frac{C\lambda}{(q^-)^{k_1+1}(k_1+1)}\|u\|_X^{q^-(k_1+1)},\\\quad \mbox{ if } \lambda>0, \beta<0. \\\\

-\frac{\lambda\lambda_1(x)}{(q^+)
^{k_1+1}(k_1+1)}\|u\|^{q^+(k_1+1)}- \frac{b}{(p^-)^{\gamma+1}(\gamma+1)}\|u\|^{p^-(\gamma+1)}_X -\frac{C\beta}{{(r^-)^{k_2+1}}(k_2+1)}\|u\|^{r^-(k_2+1)},\\\quad \mbox{ if } \lambda<0, \beta>0. \\\\

-\frac{\lambda\lambda_1(x)}{(q^+)^{k_1+1}(k_1+1)} \|u\|^{q^+(k_1+1)}-\frac{\beta\lambda_1(x)}{(r^+)^{k_2+1}(k_2+1)}\|u\|^{r^+(k_2+1)}- \frac{b}{(p^-)^{\gamma+1}(\gamma+1)}\|u\|^{p^-(\gamma+1)}_X,\\\quad \mbox{ if } \lambda<0, \beta<0.
\end{array}
\right.
\end{eqnarray*}
Thus,
\begin{eqnarray*}
\mathcal{I}_{\lambda, \beta}(u)\geq \left\{
 \begin{array}{ll}
\rho^{r^+(k_2+1)}\left(-\frac{\beta\lambda_1(x)}{(r^+)^{k_2+1}(k_2+1)}- \frac{b}{(p^-)^{\gamma+1}(\gamma+1)}\rho^{p^-(\gamma+1)-r^+(k_2+1)}\right)
-\frac{C\lambda}{(q^-)^{k_1+1}(k_1+1)}\rho^{q^-(k_1+1)},\\\quad \mbox{ if } \lambda>0, \beta<0. \\\\
\rho^{q^+(k_1+1)}\left(-\frac{\lambda\lambda_1(x)}{(q^+)^{k_1+1}(k_1+1)}- \frac{b}{(p^-)^{\gamma+1}(\gamma+1)}\rho^{p^-(\gamma+1)-q^+(k_1+1)}\right) -\frac{C\beta}{{(r^-)^{k_2+1}}(k_2+1)}\rho^{r^-(k_2+1)},\\\quad \mbox{ if } \lambda<0, \beta>0. \\\\

\rho^{r^+(k_2+1)}\left(-\frac{\beta\lambda_1(x)}{(r^+)^{k_2+1}(k_2+1)}- \frac{b}{(p^-)^{\gamma+1}(\gamma+1)}\rho^{p^-(\gamma+1)-r^+(k_2+1)}\right),\\\quad \mbox{ if } \lambda<0, \beta<0.
\end{array}
\right.
\end{eqnarray*}
Choosing
\begin{eqnarray*}
\rho\in  \left\{
 \begin{array}{ll}
\left(0,\min\left(1,\left[-\frac{(p^-)^{\gamma+1}(\gamma+1)\beta\lambda_1(x)}{2b(r^+)^{k_2+1}(k_2+1)}\right]^{p^-(\gamma+1)-r^+(k_2+1)}\right)\right),\quad \mbox{ if } \lambda>0, \beta<0 \mbox{ and } \lambda<0, \beta<0, \\\\
\left(0,\min\left(1,\left[-\frac{(p^-)^{\gamma+1}(\gamma+1)
\lambda\lambda_1(x)}{2b(q^+)^{k_1+1}(k_1+1)}\right]^{p^-(\gamma+1)-q^+(k_1+1)}\right)\right),\quad \mbox{ if } \lambda<0, \beta>0,
\end{array}
\right.
\end{eqnarray*}
we deduce, for any $u\in X$ with $\|u\|_X =\rho$, that
\begin{eqnarray}\label{rel3}
\mathcal{I}_{\lambda, \beta}(u)\geq \left\{
 \begin{array}{ll}
-\frac{\beta\lambda_1(x)}{2(r^+)^{k_2+1}(k_2+1)}\rho^{r^+(k_2+1)}
-\frac{C\lambda}{(q^-)^{k_1+1}(k_1+1)}\rho^{q^-(k_1+1)},\quad \mbox{ if } \lambda>0, \beta<0. \nonumber\\\\
-\frac{\lambda\lambda_1(x)}{2(q^+)^{k_1+1}(k_1+1)}\rho^{q^+(k_1+1)} -\frac{C\beta}{{(r^-)^{k_2+1}}(k_2+1)}\rho^{r^-(k_2+1)},\quad \mbox{ if } \lambda<0, \beta>0. \\\\

-\frac{\beta\lambda_1(x)}{2(r^+)^{k_2+1}(k_2+1)}\rho^{r^+(k_2+1)},\quad \mbox{ if } \lambda<0, \beta<0.
\end{array}
\right.
\end{eqnarray}
Now, we put
\begin{eqnarray*}
\left\{
 \begin{array}{ll}
 \lambda^*=
-\frac{\beta\lambda_1(x)(q^-)^{k_1+1}(k_1+1)}{4C(r^+)^{k_2+1}(k_2+1)\rho^{q^-(k_1+1)-r^+(k_2+1)}}
,\quad \mbox{ if } \lambda>0, \beta<0. \\\\
\beta^*=-\frac{\lambda\lambda_1(x){(r^-)^{k_2+1}}(k_2+1)}{4C(q^+)^{k_1+1}(k_1+1)
\rho^{r^-(k_2+1)-q^+(k_1+1)}} ,\quad \mbox{ if } \lambda<0, \beta>0.
\end{array}
\right.
\end{eqnarray*}
We can conclude that for any $\lambda\in(0,\lambda^*)$ (respectively $\beta \in(0,\beta^*)$) , there exists $\alpha > 0$ such that for any $u\in X$ with $\|u\|_X =\rho$

\begin{eqnarray*}
\mathcal{I}_{\lambda, \beta}(u)\geq \left\{
 \begin{array}{ll}
-\frac{\beta\lambda_1(x)}{4(r^+)^{k_2+1}(k_2+1)}\rho^{r^+(k_2+1)}:=\alpha>0,\quad \mbox{ if } \lambda\in(0,\lambda^*), \beta<0. \nonumber\\\\
-\frac{\lambda\lambda_1(x)}{4(q^+)^{k_1+1}(k_1+1)}\rho^{q^+(k_1+1)}:=\alpha>0,\quad \mbox{ if } \lambda<0, \beta\in(0,\beta^*). \\\\

-\frac{\beta\lambda_1(x)}{2(r^+)^{k_2+1}(k_2+1)}\rho^{r^+(k_2+1)}:=\alpha>0,\quad \mbox{ if } \lambda<0, \beta<0.
\end{array}
\right.
\end{eqnarray*}
We have completed the proof of Lemma \ref{lemme1K}.\qed

\begin{lem}\label{lemme1a=0}
Assume that \eqref{cond(1.4)a=0} holds. Then for every finite dimensional subspace $W\subset X$, there exists $R=R(W)>0$ such that $\mathcal{I}_{\lambda,\beta}(u) \leq0$, for all $u\in W$, with $\| u\|>R$.
\end{lem}
{\bf Proof.} Let $R=R(W)>1$, for all $u\in W$, with $\| u\|>R$, then, we have
\begin{eqnarray*}
\mathcal{I}_{\lambda, \beta}(u)  =  -\frac{b}{\gamma+1}\left(\sigma_{p(x,y)}(u)\right)^{\gamma+1} -\frac{\lambda}{k_1+1}\left(\int_\Omega\frac{1}{q(x)} |u|^{q(x)}dx \right)^{k_1+1}-\frac{\beta}{k_2+1}\left(\int_\Omega\frac{1}{r(x)} |u|^{r(x)}dx \right)^{k_2+1} \notag
\end{eqnarray*}
\begin{eqnarray*}
\mathcal{I}_{\lambda, \beta}(u)
&\leq &\left\{
 \begin{array}{ll}
-\frac{b}{\gamma+1}\left(\sigma_{p(x,y)}(u)\right)^{\gamma+1} -\frac{\lambda}{k_1+1}\frac{1}{(q^+)^{k_1+1}}\left(\int_\Omega|u|^{q(x)}dx \right)^{k_1+1}-\frac{\beta}{k_2+1}\frac{1}{(r^-)^{k_2+1}}\left(\int_\Omega |u|^{r(x)}dx \right)^{k_2+1},\\ \mbox{ if } \lambda>0, \beta<0. \\\\

-\frac{b}{\gamma+1}\left(\sigma_{p(x,y)}(u)\right)^{\gamma+1} -\frac{\lambda}{k_1+1}\frac{1}{(q^-)^{k_1+1}}\left(\int_\Omega |u|^{q(x)}dx \right)^{k_1+1}-\frac{\beta}{k_2+1}\frac{1}{(r^+)^{k_2+1}}\left(\int_\Omega |u|^{r(x)}dx \right)^{k_2+1},\\ \mbox{ if } \lambda<0, \beta>0. \\\\

-\frac{b}{\gamma+1}\left(\sigma_{p(x,y)}(u)\right)^{\gamma+1} -\frac{\lambda}{k_1+1}\frac{1}{(q^-)^{k_1+1}}\left(\int_\Omega |u|^{q(x)}dx \right)^{k_1+1}-\frac{\beta}{k_2+1}\frac{1}{(r^-)^{k_2+1}}\left(\int_\Omega |u|^{r(x)}dx \right)^{k_2+1},\\ \mbox{ if } \lambda<0, \beta<0.
\end{array}
\right.
\end{eqnarray*}
Therefore, as a consequence, all norms on the finite-dimensional space $W$ are equivalent, implying the existence of a positive constant $C_W$ such that
$$\left(\int_\Omega|u|^{q(x)}dx \right)^{k_1+1}\geq C_W \|u\|_X^{q^-(k_1+1)}\;\mbox{ and } \left(\int_\Omega|u|^{r(x)}dx \right)^{k_2+1}\geq C_W \|u\|_X^{r^-(k_2+1)}.$$
Therefore, we obtain
\begin{eqnarray*}
\mathcal{I}_{\lambda, \beta}(u)
&\leq &\left\{
 \begin{array}{ll}
-\frac{b}{(p^+)^{\gamma+1}(\gamma+1)}\|u\|^{p_-(\gamma+1)}_X -\frac{\lambda}{k_1+1}\frac{C_W}{(q^+)^{k_1+1}}
\|u\|_X^{q^-(k_1+1)}-\frac{\beta}{k_2+1}\frac{C}{(r^-)^{k_2+1}}\|u\|_X^{r^+(k_2+1)},\\ \mbox{ if } \lambda>0, \beta<0. \\\\

-\frac{b}{(p^+)^{\gamma+1}(\gamma+1)}\|u\|^{p_-(\gamma+1)}_X-\frac{\lambda}
{k_1+1}\frac{C}{(q^-)^{k_1+1}}\|u\|_X^{q^+(k_1+1)}
-\frac{\beta}{k_2+1}\frac{C_W}{(r^+)^{k_2+1}}\|u\|_X^{r^-(k_2+1)},\\ \mbox{ if } \lambda<0, \beta>0. \\\\

-\frac{b}{(p^+)^{\gamma+1}(\gamma+1)}\|u\|^{p_-(\gamma+1)}_X-\frac{\lambda}{k_1+1}
\frac{C}{(q^-)^{k_1+1}}\|u\|_X^{q^+(k_1+1)}-\frac{\beta}{k_2+1}\frac{C}{(r^-)^{k_2+1}}\|u\|_X^{r^+(k_2+1)},\\ \mbox{ if } \lambda<0, \beta<0.
\end{array}
\right.
\end{eqnarray*}

Then, it is deduced from \eqref{cond(1.4)a=0} that  $\mathcal{I}_{\lambda, \beta}(u)< 0$. Hence, the proof of Lemma \ref{lemme1a=0} is complete.\qed

\end{document}